\pdfoutput=1
\documentclass{article}
\usepackage[margin=25mm]{geometry}
\usepackage{amsmath}
\usepackage{amsfonts}
\usepackage{amssymb}
\usepackage{graphicx}
\usepackage[utf8]{inputenc}
\usepackage[T2A]{fontenc}
\usepackage[ukrainian,russian,english]{babel}
\usepackage[unicode, pdftex]{hyperref}
\usepackage{xcolor}
\usepackage[all]{xy}

\newtheorem{theorem}{Theorem}

\title{Structure of the codimesion one gradient flows with at most six singular points on the M\"obius strip}

\author{Maria Loseva, Alexandr Prishlyak, Kateryna Semenovych, Yuliia Volianiuk}

\begin{document}
\newcounter{contnumfig}
\setcounter{contnumfig}{0}

\maketitle
\begin{abstract}
   We describe all possible topological structures of Morse flows and typical one-parametric gradient bifurcation on the M\"obius strip in the case that the number of singular point of flows is at most six. To describe structures, we use the separatrix diagrams of flows. The saddle-node bifurcation is specified by selecting a separatrix in the diagram of the Morse flow befor the bifurcation and the saddle connection is specified by a separatrix, which connect two saddles on the diagram.
\end{abstract}
\section*{Introduction}

We consider gradient flows on the M\"obius strip. Since the function increases along each trajectory, the flow has no cycles and polycycles. In general position, a typical gradient flow is a Morse flow (Morse-Smale flow without closed trajectories). In typical one-parameter families of gradient flows, two types of bifurcations are possible: saddle-node and saddle connection. The vector fields at the moment of the bifurcation completely determine the topological type of the bifurcation in our case.
To classify Morse flows, we use a separatrix diagram, in which separetrices are trajectories of one-dimensional stable or unstable manifolds.  

Without loss of generality, we assume that under bifurcation (as the parameter increases), the number of singular points does not increase. The saddle-node bifurcation is defined by a separatrix, which is contructed to a point. We mark this separatrix on the diagram. A saddle connection bifurcation in the diagram corresponds to a separatrix, which conect  two saddles. 

We colar  stable separatrices in red, unstable separatrices in green and saddle connections in black.

 Reeb \cite{Reeb1946} construct topological invariants of functions oriented 2-maniofolds. It was generlized in \cite{lychak2009morse} for the case of  non-orientable two-dimensional manifolds and in   \cite{Bolsinov2004, hladysh2017topology, hladysh2019simple} for manifolds with boundary, in \cite{prishlyak2002morse} for non-compact manifolds. 

Since Morse flows as gradient flows of Morse functions, we can fix the value of functions in singular points. Then the flow determinate the topological structure of the function \cite{lychak2009morse, Smale1961}. Therefore, Morse flows structure with order of critical values determinates the structure of the functions.

Possible structures of smooth function on closed 2-manifolds was described  in \cite{bilun2023morseRP2, hladysh2019simple, hladysh2017topology,  prishlyak2002morse, prishlyak2000conjugacy,  prishlyak2007classification, lychak2009morse, prishlyak2002ms, prish2015top, prish1998sopr,  bilun2002closed,  Sharko1993, prish2002Met}, on 2-manifolds with the boundary in \cite{hladysh2016functions, hladysh2019simple, hladysh2020deformations} and on closed 3- and 4-manifolds in  \cite{prishlyak1999equivalence, prishlyak2001conjugacy, hlp2023}.

In \cite{bilun2023gradient, Kybalko2018, Oshemkov1998, Peixoto1973, prishlyak1997graphs, prishlyak2020three,  prishlyak2022topological, prishlyak2017morse,  kkp2013,  prish2002vek,  prishlyak2021flows,  prishlyak2020topology,   prishlyak2019optimal, prishlyak2022Boy, lpss2024}, 
the structures  of flows on closed 2- manifolds and 
\cite{bilun2023discrete, loseva2016topology, prishlyak2017morse, prishlyak2022topological, prishlyak2003sum, prishlyak2003topological, prishlyak1997graphs, prishlyak2019optimal} on manifolds with the boundary were investigated.
Topological properties of Morse-Smale flows on 3-manifolds was considered in \cite{pbp2023, prish1998vek,  prish2001top, Prishlyak2002beh2, prishlyak2002ms,   prish2002vek, prishlyak2005complete, prishlyak2007complete, hatamian2020heegaard, bilun2022morse, bilun2022visualization}.


The purpose of this paper is to describe all possible topological structures of the Morse flows and typical bifurcations  with no more than six singular points (a saddle-node point we consider as two points) on the M\"bious strip.



 
\section{Typical one-parameter bifurcations of gradient flows on a M\"obius strip}

Typical vector fields (flow) on compact 2-manifolds are Morse-Smale fields (flow). Morse fields (or Morse-Smale gradient-like fields) are not containing closed trajectories. They satisfy following properties:

1) it contain a finite number of singular points and its are nondegenerate;

2) there are no separatric connections between saddle points;

3) $\alpha$-limiting ($\omega$-limiting) set of each trajectory is a singular point.

In typical one-parameter field families, one of these conditions is violated. If violation of the first condition is, then we have a saddle-node bifurcation. The third condition cannot be violated for gradient fields. 

\subsection{Internal bifurcations}

According to the theory of bifurcations, there are only two typical bifurcations of gradient flows: a saddle-node and a saddle connection.

\subsubsection{Saddle-node bifurcation}
\begin{figure}[ht!]
\center{\includegraphics[width=0.95\linewidth]{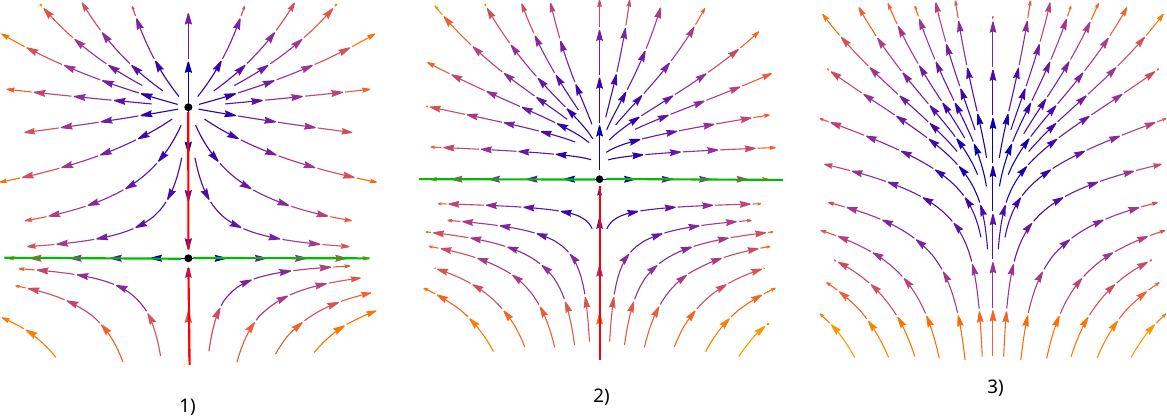}
}
\caption{SN -- saddle-node bifurcation
}
\label{bifsn}
\end{figure}

The saddle-node bifurcation, when the node is the source, is shown in Fig. \ref{bifsn}.

It can be described by the equation $V(x,y,a)=\{x,y^2+a\}$ if the node is the source and the equation $V(x,y,a)=\{-x,- y^2a\}$ if the node is a sink. Here, $a$ is a parameter. If $a<0$ we get the flow before the bifurcation, if $a>0$ we get the flow after the bifurcation, and if $a=0$ we get the flow at the moment of the bifurcation (the flow of codimensionality 1).
In order to determine the saddle-node bifurcation it is necessary to select a separatrix on the seperatrix diagram.

\subsubsection{Saddle connection}
\begin{figure}[ht!]
\center{\includegraphics[width=0.95\linewidth]{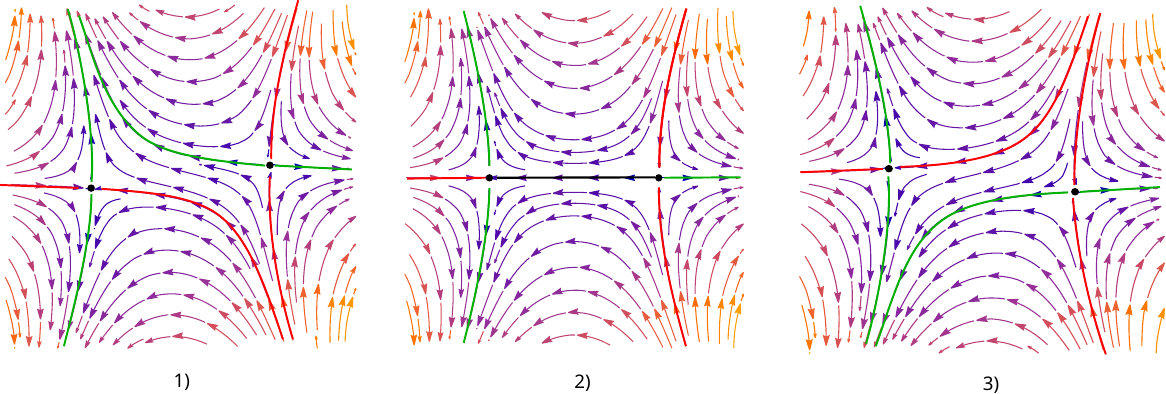}
}
\caption{SC -- saddle connection
}
\label{bifsc}
\end{figure}

The bifurcation of the saddle connection is shown in 
Fig. \ref{bifsc}. 
It can be described by the equation $V(x,y,a)=\{x^2-y^2-1,-2xy+a\}$.

\subsection{Bifurcations of singular points on the boundary}

Depending on the types of points that stick together, different options for bifurcations are possible.

\begin{figure}[ht!]
{\center{\includegraphics[width=0.9\linewidth]{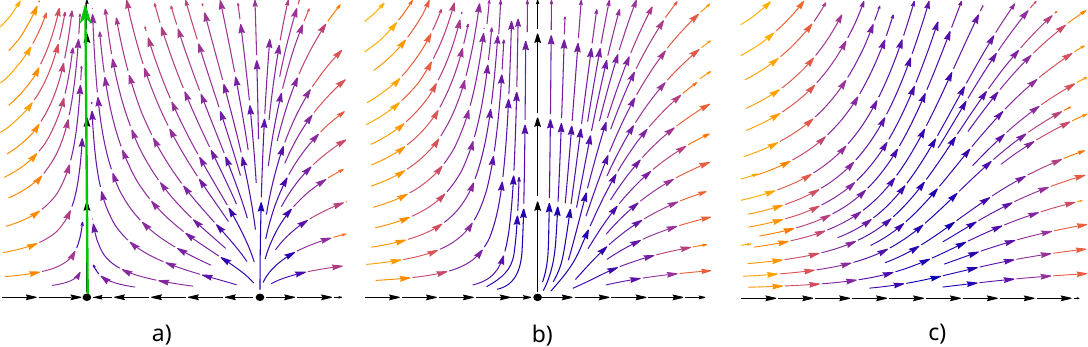}}}
\caption{ BSN -- saddle-node (source) bifurcation at the surface boundary }
\label{s-d}
\end{figure}

1) In the first case, the sourse and saddle are glueded together at a point.
In Fig. \ref{s-d} a) we show the flow before the bifurcation ($a=-1$), in Fig. \ref{s-d} b) -- flow at the moment of bifurcation ($a=0$), in Fig. \ref{s-d} c) -- flow after bifurcation ($a=1$).

\begin{figure}[ht!]
{\center{\includegraphics[width=0.9\linewidth]{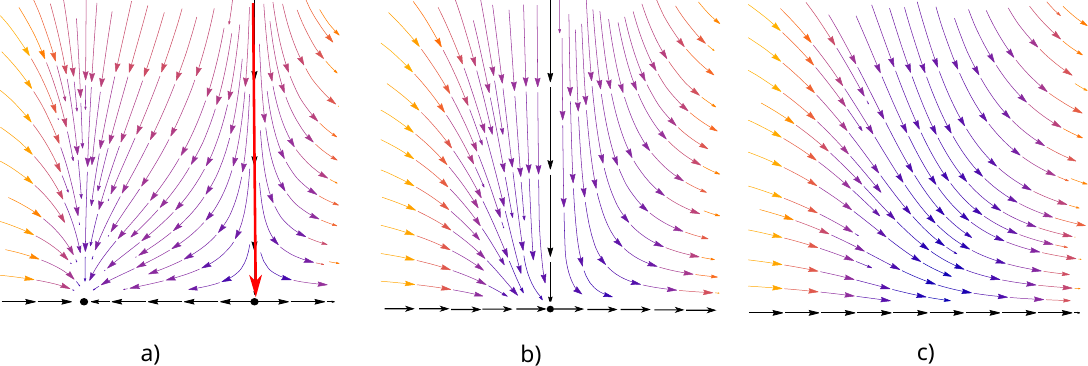}}}
\caption{ BSN -- saddle-node bifurcation (sink) at the surface boundary }
\label{s-k}
\end{figure}

2) In the second case, the saddle and the sink merge into a point, which then disappears (Fig.\ref{s-k}).

If one of the two saddle-node bifurcation points is internal, and the other lies on the boundary, then we have two types of semi-boundary saddle-node bifurcations: at the moment of bifurcation, a saddle (HS) or a node (HN).

\begin{figure}[ht!]
{\center{\includegraphics[width=0.9\linewidth]{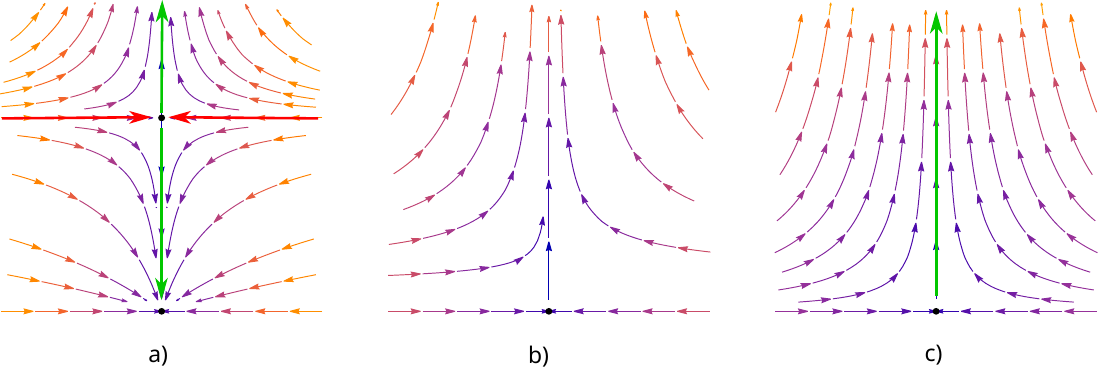}}}
\caption{ HS -- semi-boundary saddle bifurcation   }
\label{HS}
\end{figure}

\newpage

\begin{figure}[ht!]
\center{\includegraphics[width=0.9\linewidth]{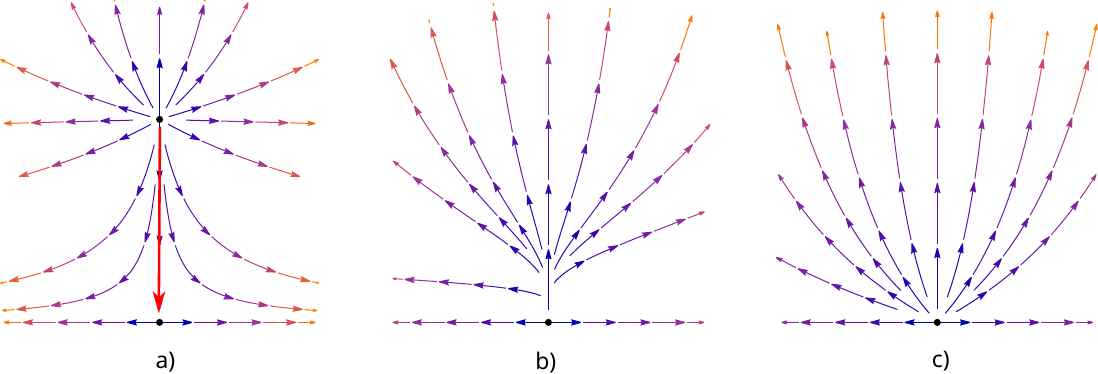}}
\label{HN}
\caption{HN -- semi-boundary node bifurcation }
\end{figure}

\begin{figure}[ht!]
\center{\includegraphics[width=0.9\linewidth]{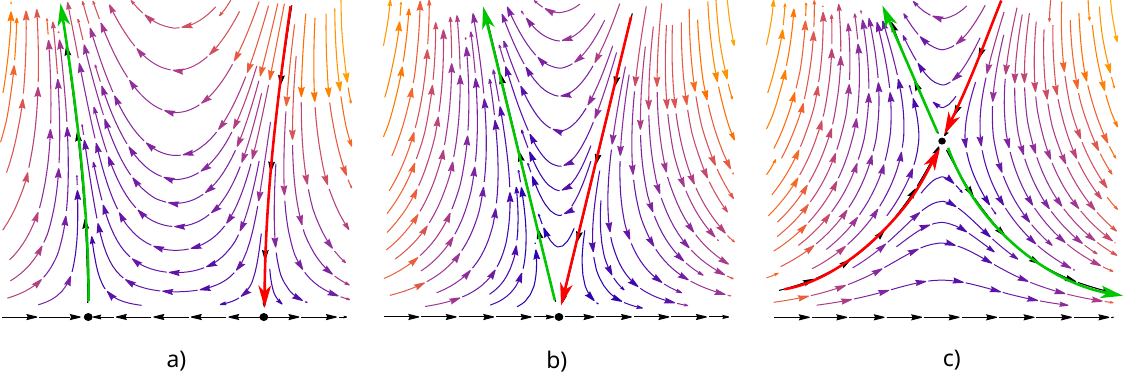}}
\caption{BDS -- double saddle bifurcation at the boundary }
\label{s-s}
\end{figure}

In addition, the following options are possible:

3) both points that stick together are saddles:
  $\{ x^2-y^2+a, -2xy \}$ (Fig. \ref{s-s}),

 On a set of flows with fixed points on the boundary and without closed trajectories, the bifurcation in a typical family is given either by the initial flow and a compressible trajectory for local bifurcations on the boundary or by the flow at the moment of bifurcation in the case of a saddle connection.

Therefore, the following types of gradient bifurcations are possible on M\"obius strip:

SN -- internal saddle node;

SC -- internal saddle connection;

BSN - boundary saddle-node;

BDS -- boundary double saddle;

HN -- semi-boundary saddle node (node);

HS -- semi-boundary saddle-node (saddle);

HSC -- semiboundary saddle connection;

BSC -- saddle connection of boundary saddle.

In the case of saddle-node bifurcations, such a bifurcation is given by a separatrix diagram to the bifurcation, on which the trajectory (separatrix) between the saddle and the node is highlighted, which is compressed to a point. To specify the bifurcation of the saddle connection, a separatrix diagram at the moment of bifurcation is sufficient.

\section{The structure of typical flows and bifurcations with no more than 4 singular points on the M\"obius strip}


To find all possible structures of Morse flows on the M\"obius strip, we use the Poincaré-Hopf theorem for doubling the vector field. Since doubling the M\"obius strip results in a Klein bottle with an Euler characteristic equal to zero, the sum of the Poincaré indices of the doubled field is also zero. Since the saddle index is -1, and the source and drain indices are 0, we have the following statement: 
the total number of Morse flow sources and sinks on the Klein bottle is equal to the number of saddles. 
Note that when doubling, internal points are doubled, but boundary points are not. Therefore, the formula for the Morse flow on the M\"obius strip is: 
\begin{equation}
2N_i+N_b=2S_i+S_b.  \label{e1}
\end{equation}

Here $N_i$ is the total number of internal sources and sinks, $N_b$ is the total number of sources and sinks on the boundary, $S_i$ is the number of internal saddles, and $S_b$ is the number of saddles on the boundary.

\begin{figure}[ht!]
\center{\includegraphics[width=0.75\linewidth]{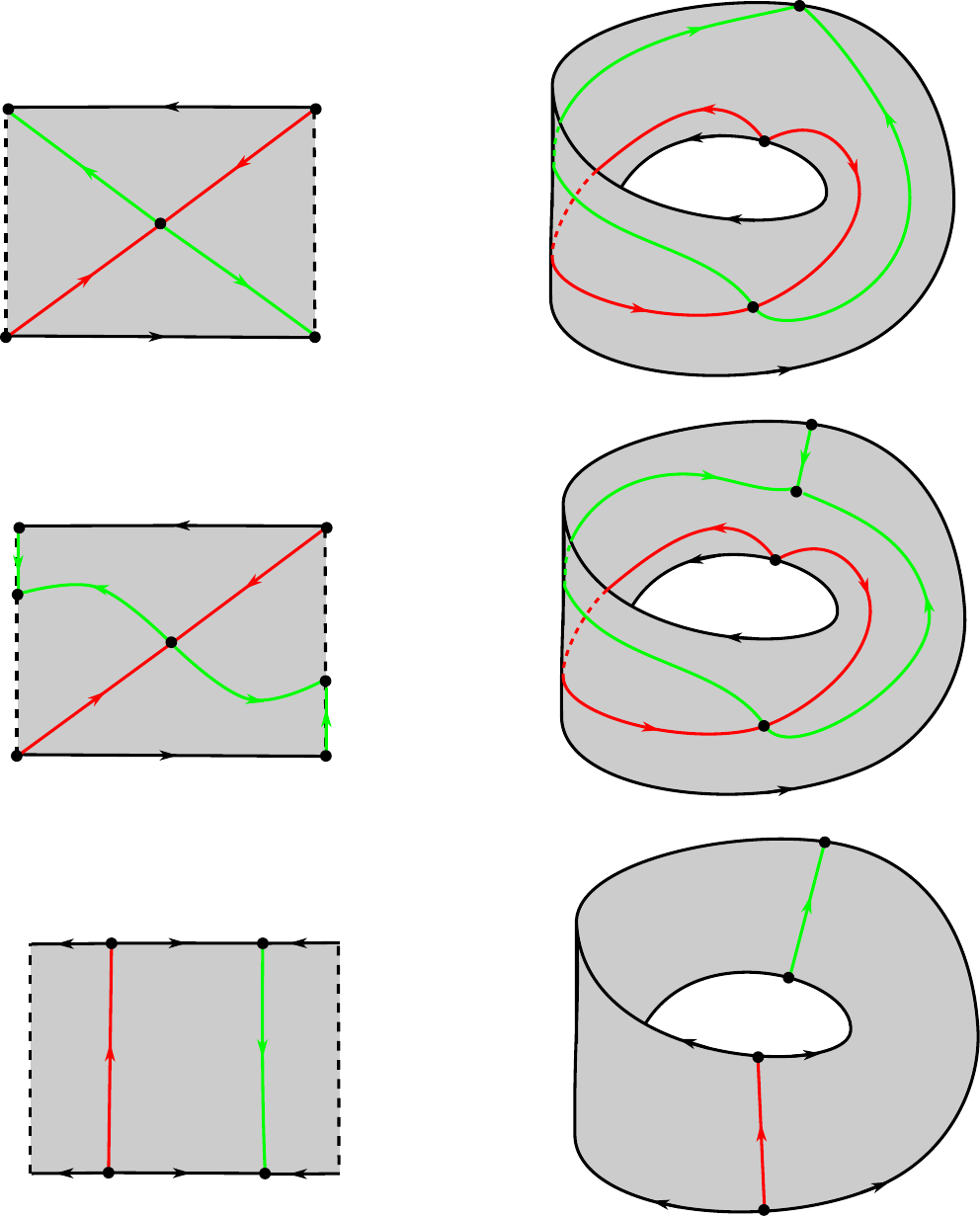}
}
\put (-200,300) {1}
\put (-200,150) {2}
\put (-200,10) {3}
\caption{Morse flows  with no more than 4 singular points}
\label{mo-ms}
\end{figure}


Еach Morse flow has a source and a sink, therefore, to fulfill the formula (1), it must contain more saddles. If the flow has three singular points, then according to (1), the only possible variant is a flow with one source and one sink at the boundary and one internal saddle. For flows with four singular points, the following options are possible: 1) internal sink and saddle, boundary source and saddle, 2) internal source and saddle, boundary sink and saddle, 3) all singular points (source, sink and two saddles) lie on the boundary.

In Fig. \ref{mo-ms}, we show all possible (with accuracy to homeomorphism) separatrix Morse flow diagrams
 with no more than 4 singular points.

Diagrams 1) and 3) are the same if we reverse orientations on the trajectories, but for diagram 2) we obtain other flow diagram. In the following figures, we will depict only one of such pair of diagrams.

For saddle-node bifurcations, it is only necessary to note how many different separatrixes and limit trajectories connecting a saddle and a node exist (with homeomorphism accuracy) on each diagram.

If the separatrix is one of the multiple edges on the separatrix diagram, then when it is pulled to a point, other multiple edges will form loops, which is not possible for gradient flows. Therefore, only those separatrices that are not one of the multiple edges should be selected.


On the diagram \ref{mo-ms}-1, all the separatrices and trajectories of the boundary are multiple edges, so it does not specify a bifurcation. In the \ref{mo-ms}-2 diagram, the only non-multiple edge is the separatrix between the saddle on the boundary and the sink. It also specifies a single bifurcation of the HN type. We get a bifurcation of the same type if we consider the reverse flow. On the \ref{mo-ms}-2 diagram, only the boundary trajectories set bifurcations -- two BSN bifurcations and one BDS.

\begin{figure}[ht!]
\center{\includegraphics[width=0.55\linewidth]{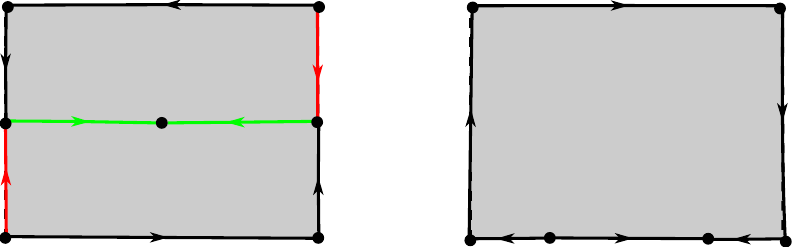}
}
\put (-200,-10) {1}
\put (-50,-10) {2}
\caption{Saddle connections in flows  with 4 singular points}
\label{mo-sc4}
\end{figure}


Saddle bifurcation is not possible for flows with three singular points, because such flows have only one saddle. For flows with two saddles, two types of saddle bifurcations are possible: HSC, if one saddle is internal and the other on the boundary, and BSC, if both saddles belong to the boundary. All possible diagrams of such flows are shown in fig.\ref{mo-sc4}. If we reverse the direction of movement in the first flow, we get a new flow, and for the second flow, the flow is equivalent to itself. 

Summarizing all of the above, we have the following:

\begin{theorem} On the Möbius strip, there exists, up to topological equivalence, a single Morse flow with three critical points.
With four critical points, there are four Morse streams and the following bifurcations: two HN bifurcations, two BSNs, one BDS, two HSCs, and one BSC.
\end{theorem}

\section{Flows and bifurcations with 5 singular points}
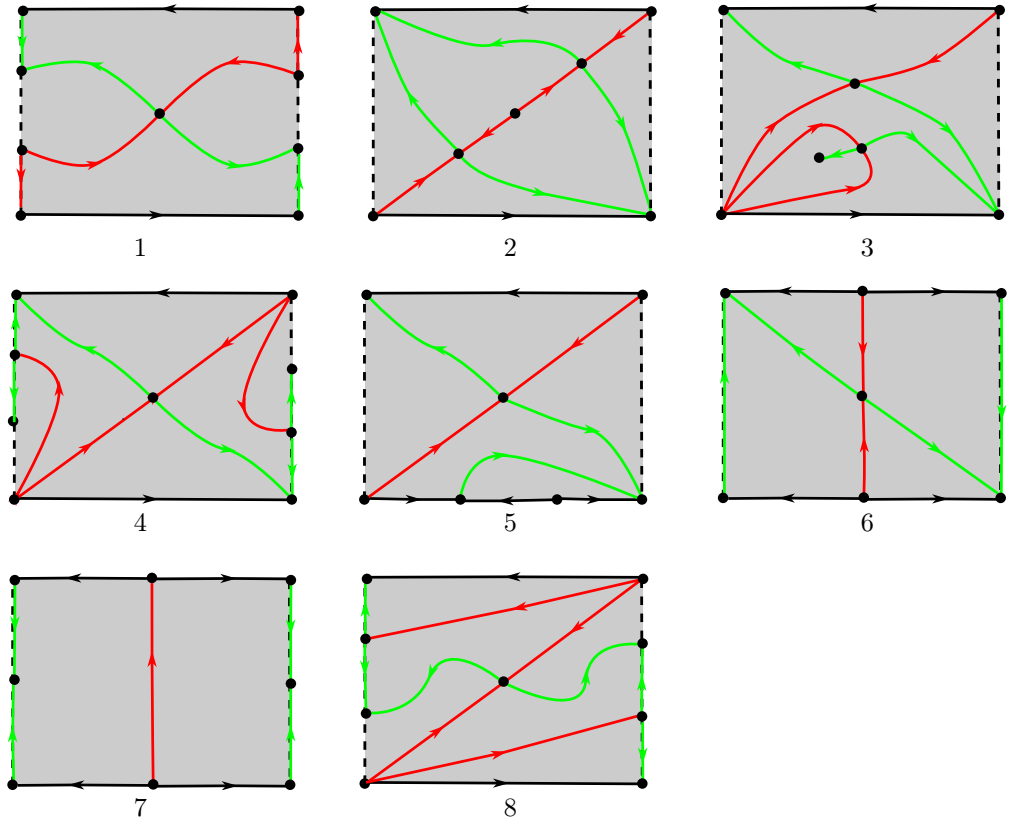
\begin{figure}[ht!] 
\center{\includegraphics[width=0.80\linewidth]{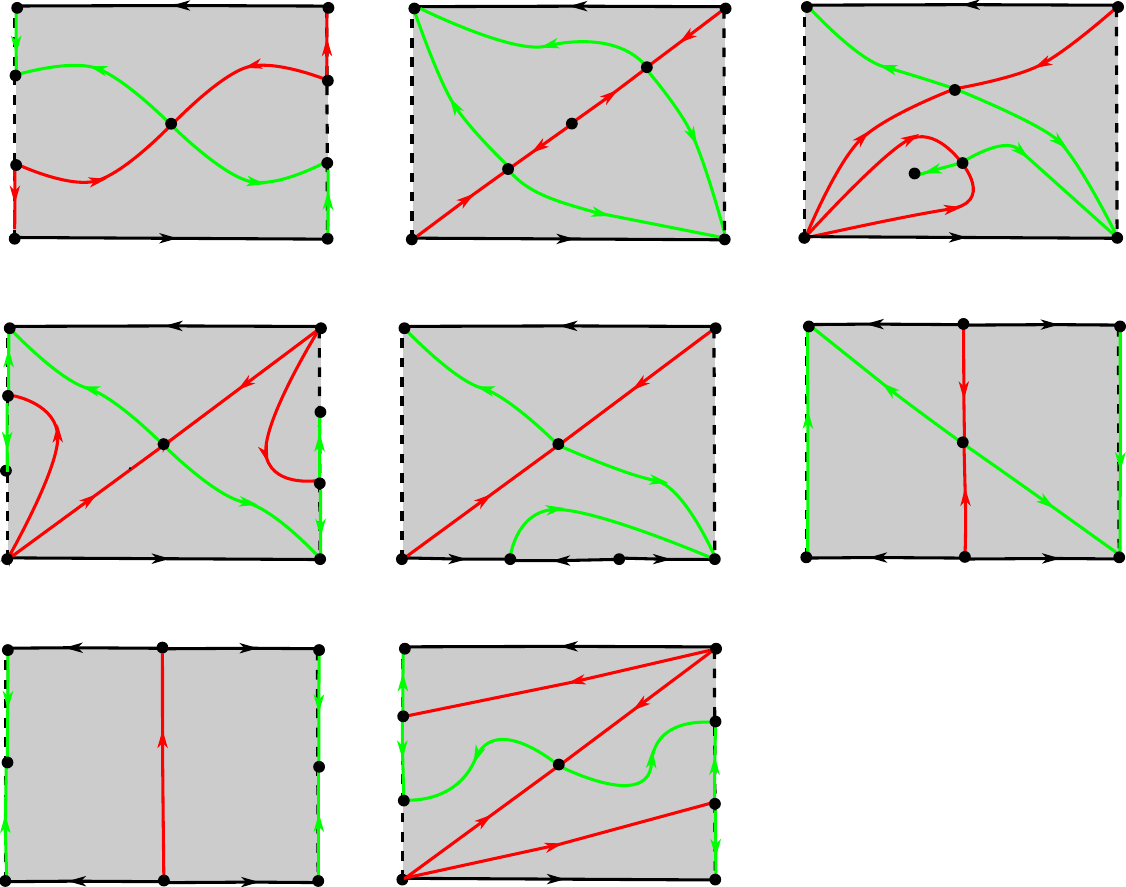}
}
\put (-330,202) {1}
\put (-190,202) {2}
\put (-55,202) {3}
\put (-330,98) {4}
\put (-190,98) {5}
\put (-55,98) {6}
\put (-330,-10) {7}
\put (-190,-10) {8}
\caption{Morse flows with 5 singular point }
\label{Mo-ms2}
\end{figure}

Next, we consider flows with five singular points (Fig. \ref{Mo-ms2}).

Since there can be only an even number of singular points on the boundary, there are either 2 or 4 of them. Let us first consider the flows with two points on the boundary. If one of these points is a saddle, and the other is not, then formula (1) is not fulfilled. Either both of these points are saddle points, or both are not saddle points. If both saddle points on the boundary are saddles, then the three interior points are the source, saddle and sink. Since in this case the separatrices can be uniquely drawn from the saddles to the source and sink (with accuracy up to homeomorphism), there is a single Morse flow, the diagram of which is shown in Fig. \ref{Mo-ms2}-1.

Let us now consider the case when both points on the boundary are not saddle points. Then one of them is a source, and the other is a sink. Internal singular points are two saddles and a node (source or sink). If three separatrices out of the node, then two of them are separatrixes of the same saddle, and therefore form a loop which is central line of M\"obius strip \ref{Mo-ms2}-8, othewise  there is another singular point inside of the loop, which is impossible. If two separators come out of the node, then the only possible flow has the diagram in Fig. \ref{Mo-ms2}-1. If the node is connected to the saddle by one separatrix, then it lies inside the loop opened by the other two separatrixes of this saddle. Two options are possible: this loop lies inside the corner adjacent to the boundary \ref{Mo-ms2}-3 or does not lie in such a corner \ref{Mo-ms2}-4.


Let's consider the case of four singular points on the boundary. Then it follows from formula (1) that one or three saddles lie on the boundary. If there is one saddle on the boundary, then the other saddle is internal, and three more non-saddle points lie on the boundary. Let, for certainty, two of them are sources, and one is a sink. There are two possibilities for red separatrixes entering the inner saddle: 1) they start at a point (Fig. \ref{Mo-ms2}-5), or 2) at different points (Fig. \ref{Mo-ms2} -6).

Let's consider the case of three saddle points on the boundary. Let the fourth point on the boundary be the source. Then the interior singular point is a sinkeds. The only possible flow has the diagram in Fig. \ref{Mo-ms2}-7. Since we have exhausted all possible options, there are no other flows with five special points on the Möbius strip.

The following saddle-node bifurcations are possible for these flows:

1)  2 HN;

2) SN, 1 HS; ($\times$2)

3) 1 SN, 1 HS; ($\times$2)

4) 1 SN, 1 HS; ($\times$2)

5) 2 BSN; ($\times$2)

6) 1 HS, 1 BSN; ($\times$2)

7) 1 HN, 1 BSN, 1 BDS; ($\times$2)

8) 1 SN, 1 HS ($\times$2).

Note that only 1) of the considered diagrams will turn into itself when the flow is reversed. Therefore, we leave it unchanged the total number of bifurcations in this case, we multiplied it  by 2 in other seven cases.

\begin{figure}[ht!]
\center{\includegraphics[width=0.80\linewidth]{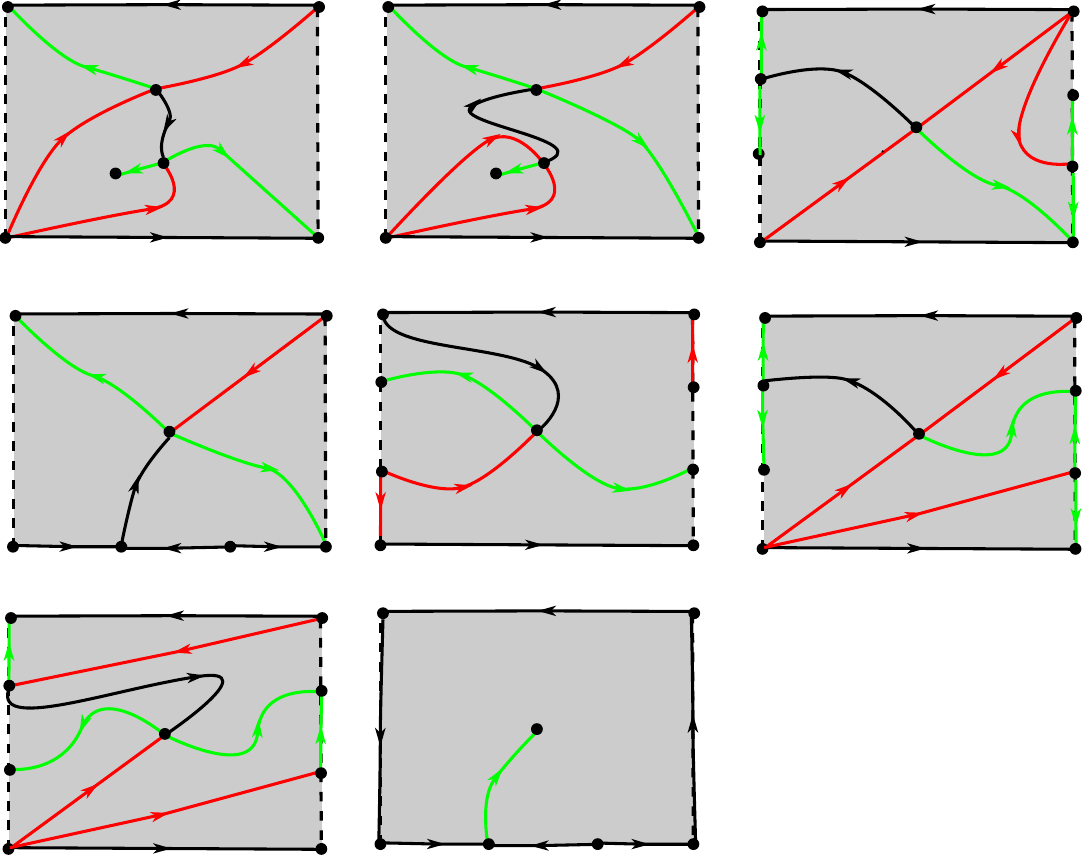}
}
\put (-330,202) {1}
\put (-190,202) {2}
\put (-55,202) {3}
\put (-330,98) {4}
\put (-190,98) {5}
\put (-55,98) {6}
\put (-330,-10) {7}
\put (-190,-10) {8}
\caption{Separatrix connections  with no more than 5 singular points}
\label{sc5}
\end{figure}
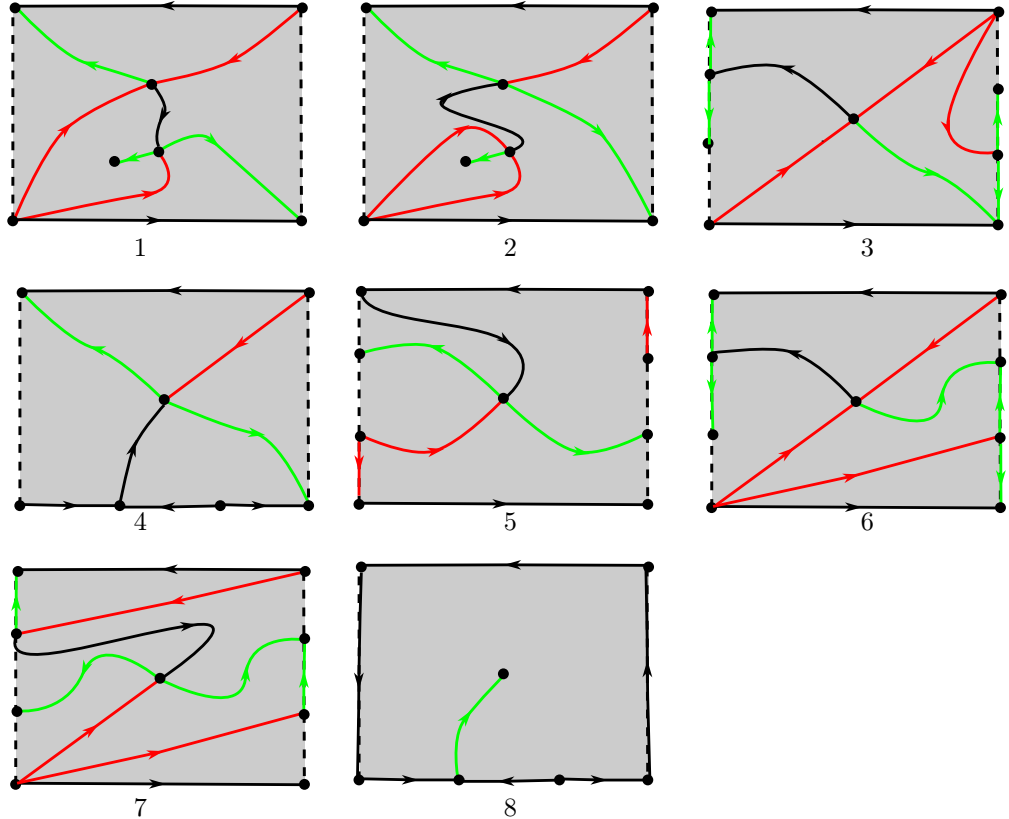

All possible separatrix connections on M\"obius stri with no more than 5 singular points are shown in fig. \ref{sc5}.

Consider flows with an internal separatrix connection (SC). In addition to the two saddles, there is another internal singular point. Let, for certainty, it be a sink. Then the singular points on the boundary are the source and the sink. Consider the case when one separatrix enters the internal sink. Diagrams of three possible flows in this case are shown in fig. \ref{sc5}.1--3. If this source includes two separators, then the possible options are shown in fig. \ref{sc5}. 6, 7. In the case of a separatrix connection between the internal and boundary saddle points (HSC), two cases are possible: 1) the boundary contains two singular points \ref{sc5}.5; 2) the boundary contains 4 singular points \ref{sc5}.4.

The only possible case of a saddle connection between points on the boundary (BSC) is shown in fig. \ref{sc5}.8. 

In all 8 cases, the inverted fields are not topologically equivalent to the original ones, so the total number of bifurcations is multiplied by 2.

\begin{theorem} On the Möbius strip, there exists, up to topological equivalence, 15 Morse flows with five critical points and following numbers of bifurcations:

10 SN bifurcations, 14  SC, 6 BSN, 2 BDS, 4 HN, 10 HS, 4HSC and 2 BSC.
\end{theorem}

\section{Flows and bifurcations with 6 singular points}

\begin{figure}[ht!]
\center{\includegraphics[width=0.96\linewidth]{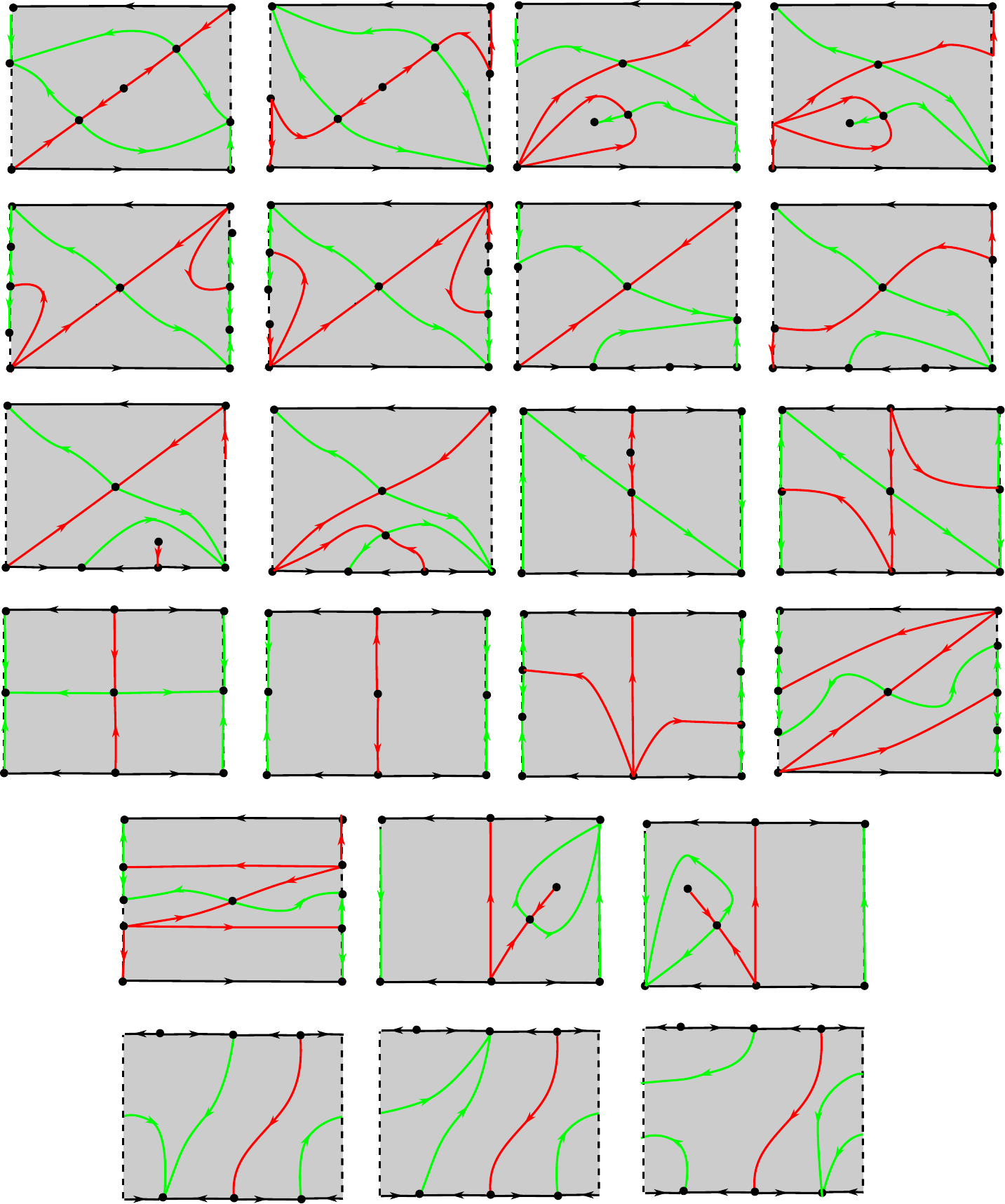}
}
\put (-400,456) {1}
\put (-285,456) {2}
\put (-170,456) {3}
\put (-55,456) {4}
\put (-400,366) {5}
\put (-285,366) {6}
\put (-170,366) {7}
\put (-55,366) {8}
\put (-400,275) {9}
\put (-285,275) {10}
\put (-170,275) {11}
\put (-55,275) {12}
\put (-400,183) {13}
\put (-285,183) {14}
\put (-170,183) {15}
\put (-55,183) {16}
\put (-350,90) {17}
\put (-235,90) {18}
\put (-120,89) {19}
\put (-350,-6) {20}
\put (-235,-6) {21}
\put (-120,-6) {22}
\caption{Morse flows with 6 singular points }
\label{ms3}
\end{figure}
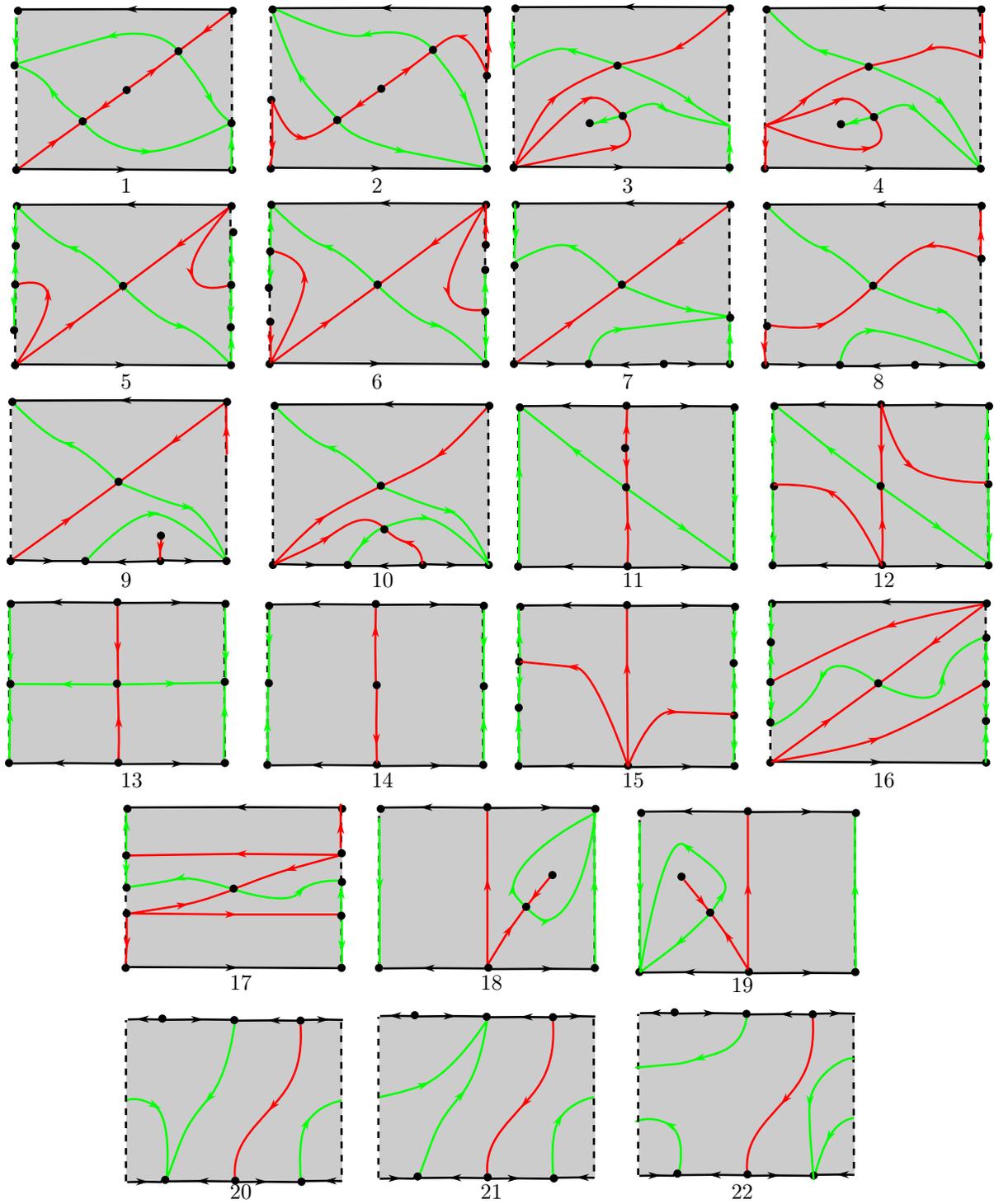

In fig. \ref{ms3} shows all possible separatrix Morse flow diagrams with 6 singular points.

The flow can have two (\ref{ms3}.1--6), four (\ref{ms3}.7--19) or six (\ref{ms3}.20--22) singular points on the boundary. 

\begin{figure}[ht!] 
\center{\includegraphics[width=0.7\linewidth]{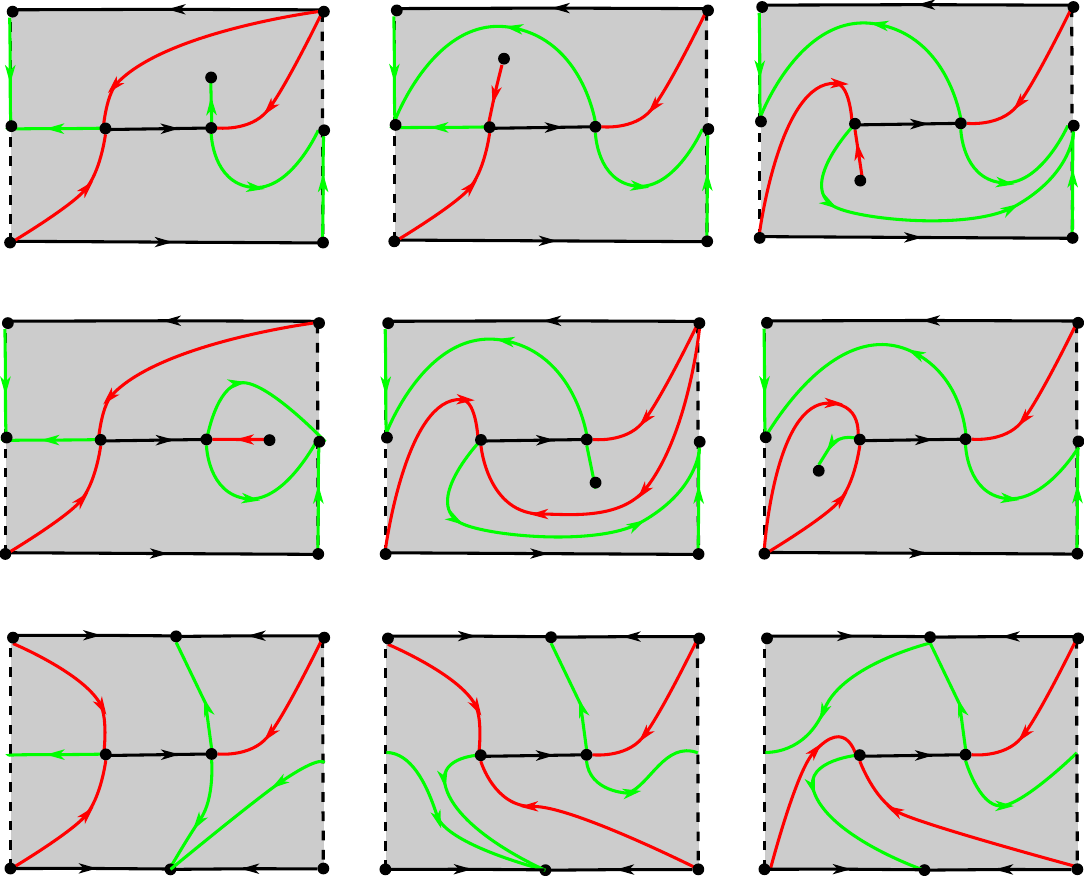}
}

\caption{SC flows with 6 singular point }
\label{Mo-sc6}
\end{figure}

\begin{figure}[ht!] 
\center{\includegraphics[width=0.7\linewidth]{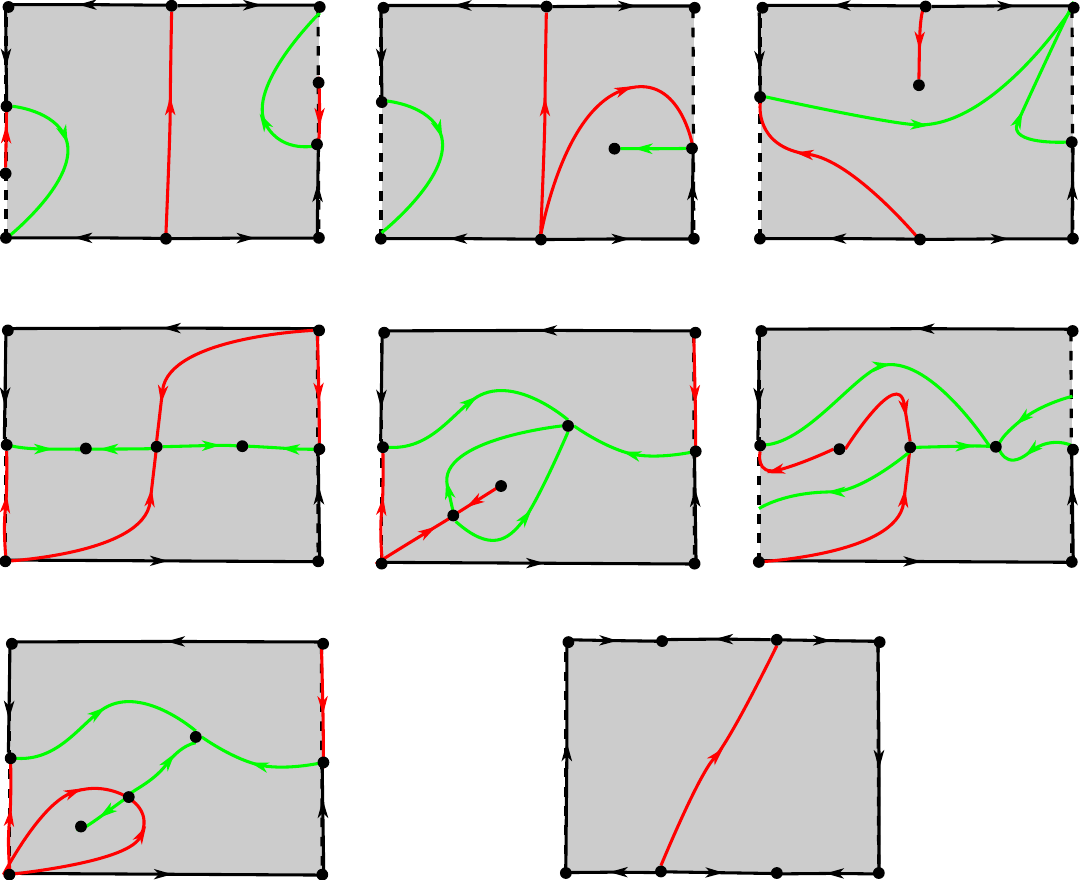}
}

\caption{HSC and BSC flows with 6 singular point }
\label{Mo-hsc6}
\end{figure}

Depending on the number of separatrixes on each of the diagrams, we get the following possible bifurcations for each of them:

1) 3 SN, 1 HN, 1 HS, (x2);

2) 2 SN, 1 HN (x2);

3) 2 SN, 1 HN, (x2);

4) 1 SN, 1 HN, 1 HS, (x2);

5) 2 SN, 1 HN,  (x2);

6) 1 SN, 1 HN, 1 HS, (x2);

7) 2 BSN,1 HN, (x2);

8) 2 BSN, 1 BDS, 1 HN, (x2);

9) 2 BSN, 1 BDS, 1 HN, (x2);

10) 4 HS;

11) 1 SN, 1 BSN, 1 HN, 1 HS, (x2);

12) 4 HS, (x2);

13) 1 BSN, 1 HN, 1 HS, (x2);

14) 1 BDS, 2 HN;

15) 1 SN, 2 BSN, 2 BDS, 1 HN, 1 HS, (x2);

16) 2 SN, 1 HN, (x2);

17) 1 SN, 1 HN, 1 HS, (x2);

18) 1 SN, 2 BSN, 1 BDS, 1 HS, (x2);

19) 1 SN, 2 BSN, 1 BDS, 1 HS, (x2);

20) 4 BSN, 1 BDS (x2);

21) 2 BSN, 2 BDS,  (x2);

22) 4 BSN, 1 BDS, (x2).

When changing the direction of the flow, diagrams 10) and 14) will not change, therefore, in the general calculation of the number of bifurcations, the above values will not change. For the rest of the charts, these values will be doubled.

Internal saddle connections are shown in Figure \ref{Mo-sc6}. In last three diagram, reverse of oreintation leads to the same diagram.
Seven HSC  and one BSC bifurcation are shown in Figure \ref{Mo-hsc6}.
In all of them,  reverse of oreintation leads to the new diagrams.

We went through all the possible options, and therefore it is fair

\begin{theorem}
The following possible structures of typical one-parameter gradient  bifurcations with 6 singular points exist on the M\"obius strip:

36SN, 15 SC, 48 BSN, 21 BDS, 30 HN, 14 HSC, 2 BSC.

\end{theorem}

\section*{Conclusion}

All possible structures of Morse flows and typical one-parameter bifurcations on M\"obius strip in which no more than six singular points are found (see Table 1). We hope that the research carried out in this paper can be extended to other surfaces and with a larger number of singular points.

\begin{table}[ht]
	\centering
		\begin{tabular} {|c|c|c|c|c|c|c|c|c|c|}
		\hline
Number of points
& 
Morse & SN & SC & BSN & BDS & HN & HS & HSC & BSC
\\ 
\hline
3 & 1 & 0 & 0 & 0 & 0 & 0 & 0 & 0 & 0 
 \\
\hline
4& 4 & 2 & 0 & 2 & 1 & 2 & 0 & 2 & 1  
 \\
\hline
5 & 15 & 10 & 14  & 6 & 2 & 4 & 10 & 4 & 2   
 \\
\hline
6& 42 & 36& 15 & 48 & 21 & 30 & 30 & 14 & 2
 \\

 \hline		

		\end{tabular}
	\caption{Number of Morse flows and bifurcations on M\"o (number of points befor bifurcation)}
	\label{tab:NF}
\end{table}


\textsc{Taras Shevchenko National University of Kyiv}

Maria Loseva \ \ \ \ \ \ \  \ \ \ \ \ \ \ \
\textit{Email:} \text{ mv.loseva@gmail.com} \   \ \ \   \ \ \ 
\textit{ Orcid ID:} \text{0000-0002-2282-206X}

Alexandr Prishlyak \ \   \ \ \  
\textit{Email address:} \text{ prishlyak@knu.ua} \ \ \ \ 
\textit{ Orcid ID:} \text{0000-0002-7164-807X}

Kateryna Semenovych \ \ 
\textit{Email:} \text{ kateryna.semenovych@knu.ua}  \ 
\textit{ Orcid ID:} \text{0000-0001-9717-1524}

Yuliia Volianiuk \ \ \ \ \ \ \ \ \ \ \ 
\textit{Email:} \text{ julia.volianiuk@knu.ua}  \   \ \ \  \ \ \ \ \
\textit{ Orcid ID:} \text{0009-0006-7203-9467}
\end{document}